\documentclass[letterpaper, 10 pt, conference]{ieeeconf}  
\usepackage{amsmath}
\usepackage[utf8]{inputenc} 
\usepackage{amsfonts}  
\usepackage{amsmath,amssymb}
\usepackage[utf8]{inputenc}
\IEEEoverridecommandlockouts   
\usepackage{tikz}
\usetikzlibrary{positioning,arrows.meta}
\usepackage{amsmath,bm}
\usepackage{placeins} 
\usepackage[caption=false,font=footnotesize]{subfig}

\title{\LARGE \bf
Gaussian Bayesian Networks for Estimating Stiff Continuous-Discrete Stochastic Systems with Ill-Conditioned Measurements
}

\author{Priyank Behera$^{1}$ and C. Robert Kenley$^{2}$
\thanks{*This work was not supported by any organization}
\thanks{$^{1}$Priyank Behera is a student in the Department of Computer Science and Department in Statistics, Purdue University, West Lafayette, U.S.A}%
\thanks{$^{2}$C. Robert Kenley  with the Edwardson School of Industrial  Engineering, Purdue  University,
        West Lafayette, IN 48907, USA
        {\tt\small kenley@purdue.edu}}%
}

\begin{document}

\maketitle
\thispagestyle{empty}
\pagestyle{empty}

\begin{abstract}

This paper introduces a Gaussian Bayesian Network-based Extended Kalman Filter (GBN-EKF) for nonlinear state estimators on stiff and ill-conditioned continuous-discrete stochastic systems, with a further analysis on systems with ill-conditioned measurements. For most nonlinear systems, the Unscented Kalman Filter (UKF) and the Cubature Kalman Filter (CKF) typically outperform the Extended Kalman Filter (EKF). But, in state estimation of stochastic systems, the EKF outperforms the CKF and UKF. This paper aims to extend the advantages of the EKF by applying a Gaussian Bayesian Network approach to the EKF (GBN-EKF), and analyzing its performance against all three filters. The GBN-EKF does not utilize any matrix inversions. This makes the GBN-EKF stable with respect to ill-conditioned matrices. Further, the GBN-EKF achieves comparable accuracy to the EKF in stiff and ill-conditioned stochastic systems, while having lower root mean squared error (RMSE) under these conditions. 

\end{abstract}

\section{Introduction}

State estimation in continuous–discrete stochastic systems is an enduring challenge in control theory and signal processing. Classical filtering techniques such as the Extended Kalman Filter (EKF) \cite{jazwinski1970, kalman1960, kalman1961}, the Unscented Kalman Filter (UKF) \cite{julier2000}, and the Cubature Kalman Filter (CKF) \cite{arasaratnam2009} have been widely used for nonlinear models. The UKF and CKF, based on higher-order approximations, typically outperform the EKF in accuracy on smooth nonlinear problems. However, in the presence of stiff stochastic systems and ill-conditioned models, this expected performance hierarchy does not hold. Recent analyses by Kulikov and Kulikova \cite{kulikov2016, kulikov2017, kulikov2018} demonstrated that stiffness in the underlying Moment Differential Equations (MDEs) can amplify numerical instability and, counterintuitively, lead to the EKF outperforming UKF and CKF in stiff regimes. Their follow-up study \cite{kulikov2017} introduced ill-conditioned measurements as a method of comparison, which will be used in this paper as well. 

While EKF gains numerical accuracy in stiff stochastic systems, it still inherits the numerical fragility of covariance propagation and inversion. In particular, when measurements themselves are ill-conditioned, matrix inversion in the update step can destabilize all conventional continuous–discrete filters. Thus, there is a clear need for estimation methods that combine the favorable stability of the EKF framework with numerical techniques that avoid ill-conditioned matrix inversion. 

To address this gap, we propose a Gaussian Bayesian Network–based EKF (GBN-EKF). The central idea is to reframe the EKF update step in terms of conditional independence structures in a Gaussian Bayesian Network \cite{howard2005influence,kenley1986influence, shachter_gaussian_1989},\cite{10.1117/12.2321776} , \cite{ZeitzMaybeck} enabling recursive state estimation without any matrix inversion. As a result, the GBN-EKF preserves the stability properties of the EKF in stiff stochastic systems while exhibiting increased robustness when the measurement model is ill-conditioned. Numerical benchmarks on Dahlquist-type and Van der Pol oscillators confirm that the GBN-EKF achieves accuracy comparable to the EKF in stiff settings, but with consistently lower accumulated RMSE under ill-conditioned measurements. This paper  demonstrates that Bayesian Network–based formulations for filtering can provide a numerically stable alternative for challenging continuous–discrete systems.

\section{Background}

We consider the standard continuous-discrete stochastic system
\begin{align}
\mathrm{d}x(t) &= f\bigl(t,x(t)\bigr)\,\mathrm{d}t + G(t)\,\mathrm{d}w(t), \qquad t>0, \label{eq:sde}\\
z_k &= h(x_k) + v_k, \qquad k\ge 1, \label{eq:meas}
\end{align}

where $x(t)\in \mathbb{R}^n$ is the state, $f:\mathbb{R}\times\mathbb{R}^n\!\to\!\mathbb{R}^n$ is a sufficiently smooth drift function, $G(t)\in\mathbb{R}^{n\times q}$, and $w(\cdot)$ is a Brownian motion with covariance $Q(t) > 0$. Measurements $z_k\in\mathbb{R}^m$ arrive at sampling instants $t_k$ (with sampling period $\delta=t_k-t_{k-1}$), through a differentiable $h:\mathbb{R}^n\!\to\!\mathbb{R}^m$ and Gaussian noise $v_k\sim\mathcal{N}(0,R_k)$.

\subsection{Continuous--Discrete Kalman filtering}
Continuous-Discrete (CD) filters propagate only the first two moments between samples by integrating the \emph{Moment Differential Equations} (MDEs) \cite{kulikov2016}:
\begin{align}
\dot{\hat x}(t) &= f\bigl(t,\hat x(t)\bigr) \label{eq:1}
\end{align}
\begin{equation}
\begin{split}
\dot P(t) &= J\bigl(t,\hat x(t)\bigr)P(t) + P(t)J^\top\bigl(t,\hat x(t)\bigr) \\
          &\quad + G(t)Q(t)G^\top(t) \label{eq:2}
\end{split}
\end{equation}
where $J(t,\hat x)=\partial f/\partial x$ evaluated on the predicted mean.
\begin{align}\qquad \hat x(t_{k-1}) = \hat x_{k-1|k-1} \label{eq:3}, \\
P(t_{k-1})=P_{k-1|k-1} \label{eq:4}
\end{align}

Integration of these equations over $[t_{k-1},t_k]$ yields the predicted state and covariance:
\begin{equation}
\hat x_{k|k-1} = \hat x(t_k), 
\qquad 
P_{k|k-1} = P(t_k).
\end{equation}

At $t_k$ the standard EKF measurement update is applied:
\begin{align}
S_{k} &= R_k + H_k P_{k|k-1} H_k^\top, \quad K_k = P_{k|k-1} H_k^\top S_{k}^{-1}, \\
\hat x_{k|k} &= \hat x_{k|k-1} + K_k\bigl(z_k - h(\hat x_{k|k-1})\bigr), \\
P_{k|k} &= P_{k|k-1} - K_k H_k P_{k|k-1},
\end{align}
with $H_k=\partial h/\partial x$ at $\hat x_{k|k-1}$. These equations are the backbone for CD-EKF, while CD-UKF and CD-CKF use sigma/cubature constructions but still rely on MDE-type propagation. Algorithms for CD-UKF and CD-CKF can be found in \cite{kulikov2018}.

\subsection{Accumulated Root Mean Squared Error}
Our main mode of comparison is the accumulated root mean squared error (ARMSE) \cite{kulikov2016}. We evaluate estimators over a uniform grid $t_0:t_K$ using the accumulated root mean squared error (ARMSE) across $L$ Monte-Carlo runs:
\begin{equation}
\mathrm{ARMSE} \;=\;
\left(
\frac{1}{L K}\sum_{\ell=1}^{L}\sum_{k=1}^{K}
\bigl\|x^{\mathrm{ref},\,\ell}(t_k) - \hat x^{\,\ell}_{k|k}\bigr\|_2^2
\right)^{\!1/2}
\end{equation}
with identical measurement noise sequences across filters for fair comparison.

\section{Gaussian Bayesian Network-based Extended Kalman Filter}

\subsection{Gaussian Bayesian Networks}
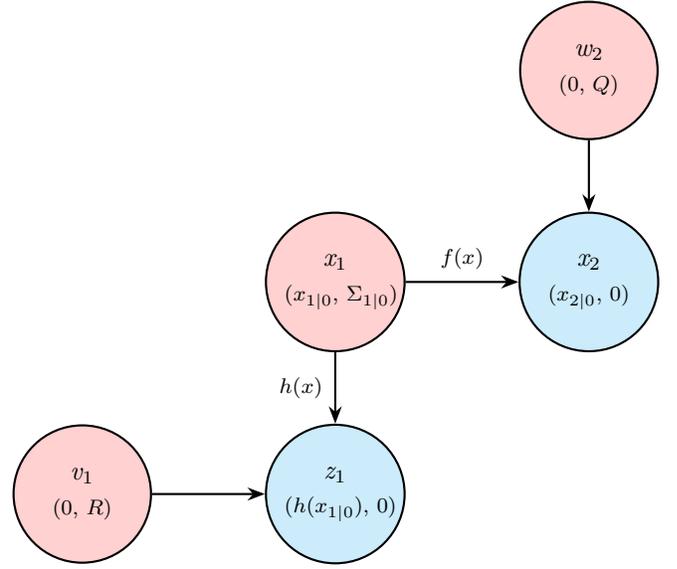
\begin{figure}[t]
    \centering
    \begin{tikzpicture}[%
        node distance=0.95cm and 1.5cm,
        >=Stealth,
        every node/.style={
            align=center,
            text width=1.35cm,
            font=\normalsize
        },
        bigc/.style={circle, draw=black, thick, fill=red!18, minimum size=1.7cm},
        bigo/.style={circle, draw=black, thick, fill=cyan!18, minimum size=1.7cm}
    ]
    \node[bigc] (x1)  {$\mathit{x}_1$\\[0.1em] {\footnotesize $({x}_{1|0},\,\Sigma_{1|0})$}};
    \node[bigo, below=of x1] (z1) {$\mathit{z}_1$\\[0.1em] {\footnotesize $(h({x}_{1|0}),\,0)$}};
    \node[bigo, right=of x1] (x2) {$\mathit{x}_2$\\[0.1em] {\footnotesize $({x}_{2|0},\,0)$}};
    \node[bigc, above=of x2] (w2) {$\mathit{w}_2$\\[0.1em] {\footnotesize $(0,\,Q)$}};
    \node[bigc, left=of z1] (v1) {$\mathit{v}_1$\\[0.1em] {\footnotesize $(0,\,R)$}};
    \draw[->, thick] (x1) -- node[above, yshift=1pt] {\footnotesize$f(x)$} (x2);
    \draw[->, thick] (x1) -- node[left, xshift=10pt, yshift=0pt] {\footnotesize$h(x)$} (z1);
    \draw[->, thick] (v1) -- (z1);
    \draw[->, thick] (w2) -- (x2);
    \end{tikzpicture}
    \caption{Gaussian Bayesian Network for discrete-time filtering. Each of these nodes is assumed to be normally distributed, taking the form $x \sim \mathcal{N}(\mu,\, \Sigma)$. The blue nodes are deterministic, and the pink nodes are stochastic
}
    \label{fig:gaussian_id}
\end{figure}

A Gaussian Bayesian Network \cite{ howard2005influence} is a probabilistic graphical model that provides a framework for inference. In this Bayesian network, each node represents a Gaussian random variable, while the edges (arcs) represent conditional linear dependencies as shown in Fig.~\ref{fig:gaussian_id}. This Bayesian network allows for probabilistic inference through arc reversals and node removals \cite{shachter_gaussian_1989, kenley1986influence}. A multivariate normal random vector $\mathbf{x} \sim \mathcal{N}(\boldsymbol{\mu}, \mathbf{\Sigma})$ can be equivalently represented in Gaussian influence diagram form by decomposing the covariance matrix $\mathbf{\Sigma}$ into (i) a strictly upper triangular matrix $\mathbf{B}$ of arc (regression) coefficients and (ii) a vector $\mathbf{V}$ of conditional variances:
\begin{align}
    x_j &= \sum_{k=1}^{j-1} B_{k j} x_k + \epsilon_j, \quad \epsilon_j \sim \mathcal{N}(0, v_j), \quad j=1, \ldots, n \\
    \mathbf{x} &\sim \mathcal{N}(\boldsymbol{\mu}, \mathbf{\Sigma}) \iff (\boldsymbol{\mu}, \mathbf{B}, \mathbf{V})
\end{align}
Here, $B_{k j}$ and $V_j$ are computed recursively from the entries of $\mathbf{\Sigma}$:
\begin{align}
    B_{k j} &= \left[\mathbf{P}_{1:(j-1),1:(j-1)} \, \mathbf{\Sigma}_{1:(j-1),j}\right]_k \\
    V_j &= \Sigma_{j j} - \mathbf{\Sigma}_{j,1:(j-1)} \mathbf{B}_{1:(j-1),j}
\end{align}
where $\mathbf{P}$ is the inverse (or generalized inverse) of $\mathbf{\Sigma}$. This conversion allows us to perform operations such as arc reversal and arc removal, using $\mathbf{B}$ and $\mathbf{V}$. 

\subsection{Time Update}

Similar to CD-EKF, the time update in CD-GBN-EKF use the same MDEs described in \eqref{eq:1}, \eqref{eq:2}, \eqref{eq:3}, \eqref{eq:4}. Typically, the time-update, or the prediction step, is similar across most continuous-discrete filters. 

\subsection{Measurement Update}

Given the predicted state $x_{k|k-1} \in \mathbb{R}^n$ and covariance $P_{k|k-1} \in \mathbb{R}^{n\times n}$,
a nonlinear measurement model $z_k = h(x_k) + v_k$, $v_k \sim \mathcal{N}(0,R)$ with $R \in \mathbb{R}^{m\times m}$,
the measurement update in CD-GBN-EKF proceeds as follows.

\begin{enumerate}
    \item \textbf{Linearization at the prior mean:}
    \begin{align}
        H_k \;\triangleq\; \left.\frac{\partial h(x)}{\partial x}\right|_{x=x_{k|k-1}} &&
        h_{k|k-1} \;\triangleq\; h(x_{k|k-1})
    \end{align}
    We define $H_k$ to be the Jacobian evaluated at the predicted state, and the $h_k$ to be the predicted measurement at the $x_{k|k-1}$.

    \item \textbf{Convert prior to Bayesian network parameters:}
    \begin{align}
        (\bm{x}_{k|k-1}, \bm{\Sigma}_{k|k-1}) \longrightarrow (\bm{x}_{k|k-1}, \bm{B}_{k|k-1}, \bm{V}_{k|k-1})
    \end{align}

    \item \textbf{Measurement Update:}
    In the Bayesian network formulation, the measurement update is performed by introducing measurement evidence and adjusting the diagram’s structure to incorporate this new information. This involves three key steps:
    \begin{enumerate}
        \item \textbf{Augmenting the System:}
    \\
    The state and measurement nodes are stacked to form an augmented mean vector and block matrix:
    \begin{align}
    \mathbf{x}_{\text{aug}} &=
    \begin{bmatrix}
        \mathbf{x}_{k|k-1} \\
        \mathbf{h}_{k|k-1}
    \end{bmatrix}, \\
    \mathbf{V}_{\text{aug}} &=
    \begin{bmatrix}
        \mathbf{V}_{k|k-1} \\
        \mathrm{diag}(\mathbf{R})
    \end{bmatrix}, \\
    \mathbf{B}_{\text{aug}} &=
    \begin{bmatrix}
        \mathbf{B}_{k|k-1}  & \mathbf{H}_k^\top \\
        \mathbf{0} & \mathbf{0}
    \end{bmatrix}
\end{align}

    where $\mathbf{R}$ is the measurement noise covariance.

    \item \textbf{Evidence Entry:}
    \\
    The measurement $\mathbf{z}$ is entered as evidence. This step utilizes Bayes' rule to update the state by conditioning the joint Gaussian on the observed value. This is performed by the 'evidence' operation
    \begin{equation}
        \texttt{evidence}\left(\mathbf{x}_{\text{aug}}, \bm{B_{\text{aug}}}, \bm{V_{\text{aug}}}, \mathbf{z}, n_0, n_1, n_2, \Delta \mathbf{x}\right)
    \end{equation}
    where $n_0, n_1, n_2$ index the state, measurement, and successor node dimensions, respectively.

    \item \textbf{Arc Reversal:}
    \\
    Evidence entry is implemented via a series of \textit{arc reversals}, which invert the direction of dependencies as needed to ensure the state nodes are conditioned on the measurement. At each step, the reversal operation updates the coefficients and variances according to:
    \begin{align}
        B'_{ji} &= \frac{1}{B_{ij}}, \\
        V'_j &= V_j + \frac{V_i}{B_{ij}^2}
    \end{align}
    for the reversed arc from node $i$ to $j$
\end{enumerate}

    \item \textbf{Convert back to covariance}
    \begin{align}
    (\bm{x}_{k|k}, \bm{B}_{k|k}, \bm{V}_{k|k}) \longrightarrow (\bm{x}_{k|k}, \bm{\Sigma}_{k|k})
    \end{align}
    $x_{k|k}$ and $\Sigma_{k|k}$ are the final estimation of the state given the measurement evidence. 
\end{enumerate}

Ill-conditioned measurement functions make conventional covariance-form updates sensitive to inversion of the innovation covariance $S$. Through this Gaussian Bayesian Network formulation, we circumvented the necessity of inversion. Further, roundoff error through inversion can lead to covariance matrices losing positive semi-definiteness, while CD-GBN-EKF ensures positive semi-definiteness.  

\section{Stability Analysis}

To illustrate the differences between conventional CD-EKF/CD-GBN-EKF and CD-UKF, we apply the described filtering methods to the Dahlquist test equation \cite{wanner1996solving} by setting  $f(t,x(t)) \equiv f(x(t))$ and $G(t) \equiv 0,$. Since we want observe filter performance across stiffness, the driving noise is effectively 0. Applying this the MDES \eqref{eq:1}, \eqref{eq:2} for the CD-EKF reduce to
\begin{align}
\hat{x}'(t) &= f\bigl(\hat{x}(t)\bigr), \label{eq:cdekf_x} \\
P'(t) &= 2\,\partial_x f\bigl(\hat{x}(t)\bigr)\,P(t), \label{eq:cdekf_p}
\end{align}
where $\partial_x f(\hat{x}(t))$ denotes the Jacobian of $f(x)$ evaluated at the mean state.

For the CD-UKF, the corresponding equations are
\begin{align}
\hat{x}'(t) &= f\bigl(\hat{x}(t)\bigr) 
+ \partial^2_{xx} f\bigl(\hat{x}(t)\bigr)\,P(t) + \text{HOT}, \label{eq:cdukf_x}\\
P'(t) &= 2\,\partial_x f\bigl(\hat{x}(t)\bigr)\,P(t) + \text{HOT}, \label{eq:cdukf_p}
\end{align}
where $\partial^2_{xx} f(\hat{x}(t))$ is the second derivative of $f(x)$ and HOT denotes higher-order terms omitted in the expansion. 

Equations~\eqref{eq:cdekf_x}–\eqref{eq:cdukf_p} reveal that both filters evolve the covariance in the same way up to second-order accuracy, but the CD-UKF includes an extra correction in the state mean. When $\partial^2_{xx} f(\hat{x}(t)) \neq 0$, this term can significantly affect performance. In stable regimes ($\partial_x f(\hat{x}(t)) \leq 0$), the covariance remains bounded, allowing the correction to improve accuracy, so the CD-UKF typically outperforms the CD-EKF. However, in unstable regimes ($\partial_x f(\hat{x}(t)) > 0$) the covariance grows rapidly, leading to unreliable computation of the correction term. In such cases, the CD-EKF is more robust. A similar argument holds for CD-CKF, which uses square root propagation. For the CD-EKF, the time-update equations are
\begin{align}
\hat{x}'(t) &= f\bigl(\hat{x}(t)\bigr), \label{eq:sr_cdekf_x}\\
S'(t) &= \partial_x f\bigl(\hat{x}(t)\bigr)\,S(t), \label{eq:sr_cdekf_s}
\end{align}
where $S(t)$ is the square-root of the covariance.

The CD-CKF equations are
\begin{align}
\hat{x}'(t) &= f\bigl(\hat{x}(t)\bigr) 
+ \tfrac{1}{2}\,\partial^2_{xx} f\bigl(\hat{x}(t)\bigr)\,S^2(t) + \text{HOT}, \label{eq:cdckf_x}\\
S'(t) &= \partial_x f\bigl(\hat{x}(t)\bigr)\,S(t) + \text{HOT}. \label{eq:cdckf_s}
\end{align}
Similar to the discussion above, the additional second-order correction in~\eqref{eq:cdckf_x} enhances accuracy under stable conditions but can be detrimental when instability causes rapid covariance growth. 

CD-GBN-EKF follows the same set of MDEs as the CD-EKF, and therefore its nominal behavior under well-conditioned measurement scenarios is essentially equivalent. The key advantage of the CD-GBN-EKF emerges in the presence of ill-conditioned measurement models. In such settings, the explicit inversion of the innovation covariance, required in the conventional CD-EKF, can lead to significant numerical instability and a subsequent degradation in estimation accuracy. By formulating the update in the information domain, the CD-GBN-EKF circumvents this inversion step, thereby maintaining numerical robustness. As a result, the filter demonstrates a substantial reduction in ARMSE relative to the CD-EKF in ill-conditioned measurement scenarios. This improvement is particularly relevant for stiff stochastic systems, where small perturbations in the measurement covariance can otherwise propagate into large estimation errors. In the first section, we first analyze filter performance under well-conditioned stiff stochastic systems. Then, we provide numerical evidence highlighting the superior performance of the CD-GBN-EKF in ill-conditioned measurement models.

\section{Numerical Study}

To ensure a consistent basis for evaluating filter performance across all test cases, we employed a uniform numerical integration scheme. Specifically, Python's \verb|solve_ivp| routine was used with parameters \texttt{rtol}$=10^{-12}$, \texttt{atol}$=10^{-12}$, \texttt{MaxStep}$=10^{-1}$, and \texttt{method="Radau"}, which is well suited for stiff differential equations. For comparative assessment, we report the ARMSE of the CD-CKF, CD-UKF, CD-EKF, and CD-GBN-EKF in both the stiff Dahlquist \cite{dahlquist1963} and stiff Van der Pol \cite{vanderpol1926} problems with well-conditioned measurements. In the case of the stiff Van der Pol system with ill-conditioned measurements, only the CD-EKF and CD-GBN-EKF were considered, as the remaining filters required prohibitively long runtimes (on the order of several hours) to complete a single trial. 

\subsection{Dahlquist-type SDE}
\begin{figure}
    \centering
    \includegraphics[width=0.42\textwidth]{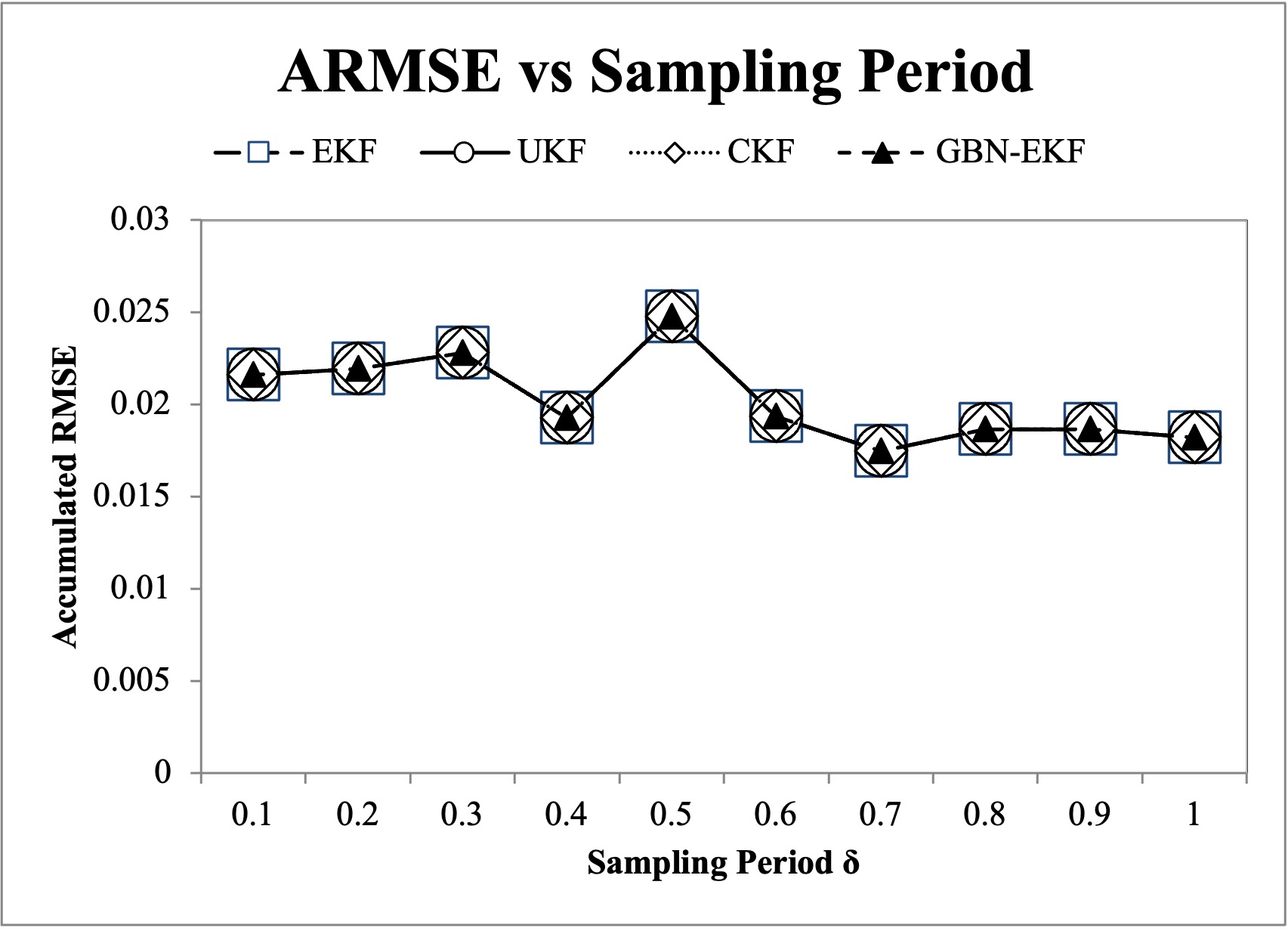}
    \caption{
        This figure plots ARMSE vs $\delta$ for the Dahlquist SDE with $j = 1$
    }
    \label{fig:Plot1}
\end{figure}

\begin{figure}
    \centering
    \includegraphics[width=0.42\textwidth]{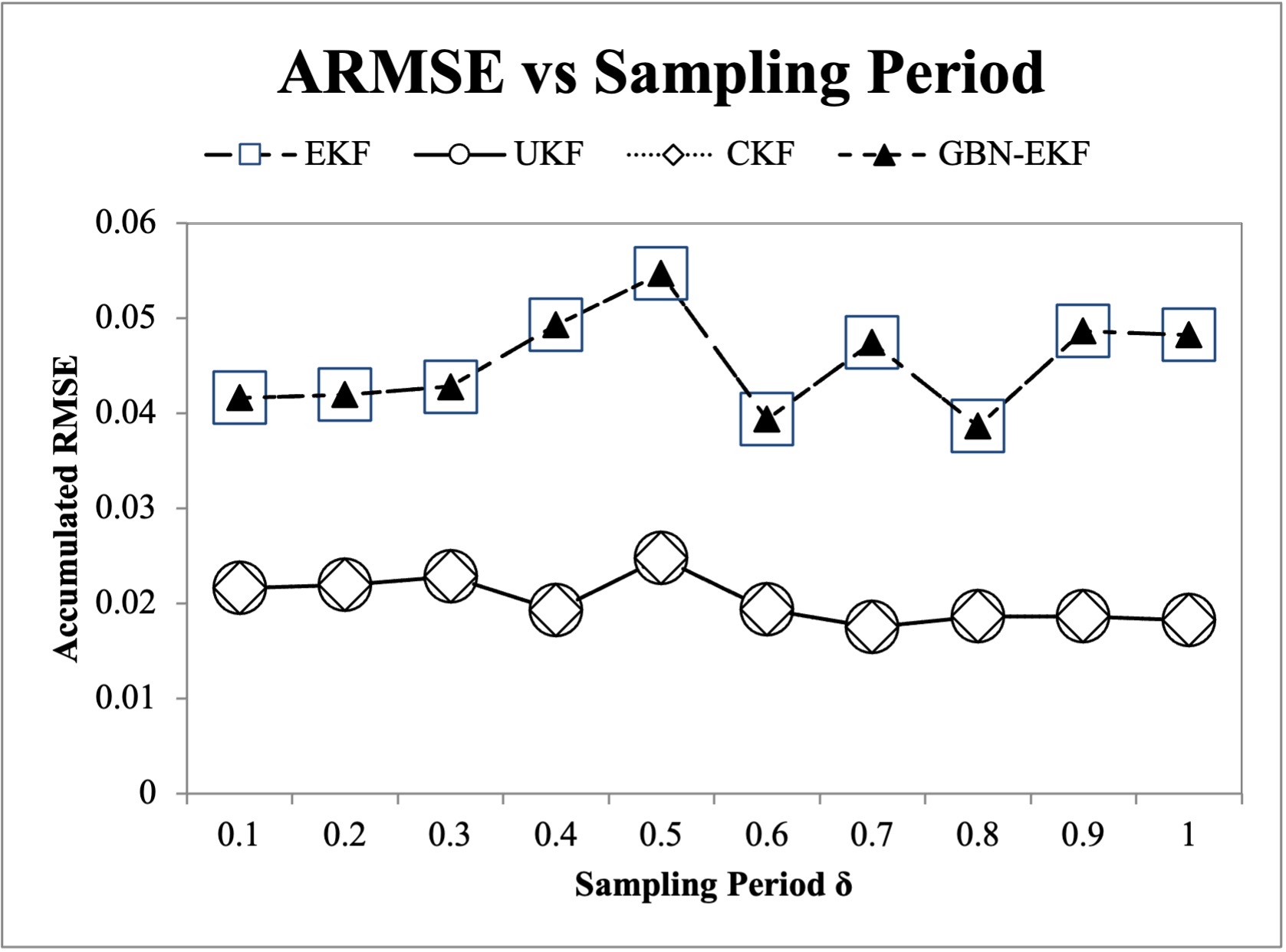}
    \caption{
        This figure plots ARMSE vs $\delta$ for the Dahlquist SDE with $j = 3$
    }
    \label{fig:Plot2}
\end{figure}

First, we estimate scalar SDE's of the Dahlquist type with the following differential equation:
\begin{align}
dx(t) \;=\; \mu x^j(t)\,dt \;+\; dw(t), 
\quad j = 1,2,3\;\; t \in [0,4].
\end{align}
and measurement equation:
\begin{align}
z_k = x(t_k) + v_k
\end{align}
Where $\mu$ is the stiffness, $w(t)$ is a Brownian motion with zero mean and unit covariance, and $j$ sets the nonlinearity of the model. The measurement noise is normally distributed with $v_k \sim \mathcal{N}(0,\,0.04)$. For the stiff case, we set $\mu=-10^{4}$, initial state $x(0) = 1$ and initial covariance $P = 10^{-2}$. Then, we plot the values for varying $\delta$, from $0.1, 0.2, ..., 1$
First, we take $j = 1$. The second-order derivative of the drift function in this type of SDE is 0, meaning that we expect the performance of all the filters to be the same. As shown in Figure \ref{fig:Plot1}, the ARMSE is the same for all filters. Next, we compute ARMSE with $j = 3$. Since the second-order derivative $\partial_{xx}^2f\bigl(\hat{x}(t)\bigr) \leq 0$, we expect CD-CKF and CD-UKF to perform better than CD-EKF and CD-GBN-EKF. Specifically, $\partial_{xx}^2f\bigl(\hat{x}(t)\bigr) = 6\mu\hat{x}(t)$, and since we set $\mu = -10^4$, we expect lower covariance in CD-CKF and CD-UKF, which is shown in Fig \ref{fig:Plot2}. 

\subsection{Van der Pol oscillator}
For this section, we consider the Van der Pol oscillator, which is modeled by the following SDE. \begin{align}
d\begin{bmatrix} x_1(t) \\ x_2(t) \end{bmatrix}
&=
\begin{bmatrix}
x_2(t) \\
\mu\bigl((1 - x_1^2(t))\,x_2(t) - x_1(t)\bigr)
\end{bmatrix} dt \notag\\
&\quad+
\begin{bmatrix}
0 & 0 \\
0 & 1
\end{bmatrix}
\, dw(t).
\end{align}

We set the initial state vector $x(0) = [1,0]^\top$ and initial covariance to be the diagonal matrix $P_0=diag\{0.05, 0.05\}$. We set $\mu = 10^4$ to make the SDE stiff. Since the second-order derivative $\partial_{xx}^2f > 0$, the state estimation because quite poor for CD-UKF and CD-CKF. We first present the results for the well-conditioned oscillator followed by the ill-conditioned one. 

\subsubsection{Well-Conditioned Measurement Model}
\begin{figure}
    \centering
    \includegraphics[width=0.42\textwidth]{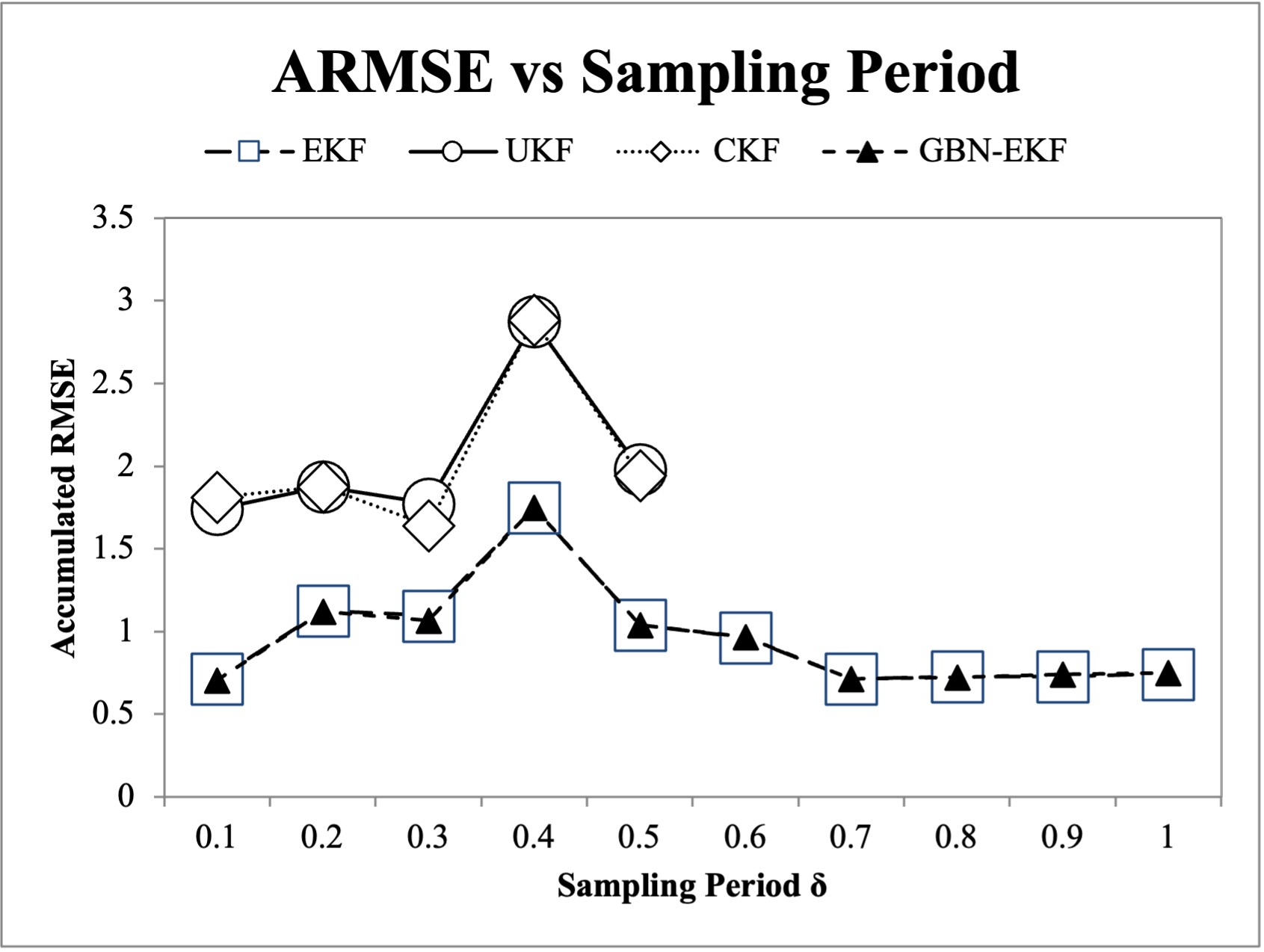}
    \caption{
        This figure plots ARMSE vs $\delta$ for the Van der Pol oscillator SDE with well-conditioned measurements
    }
    \label{fig:Plot3}
\end{figure}

For the well-conditioned measurement model, we assume full observability of the state, resulting in the following measurement equation: $z_k = x_1(t_k )+ x_2(t_k )+ v_k$, where $v_k \sim \mathcal{N}(0,\,0.04)$ is the process noise. We find that CD-EKF and CD-GBN-EKF outperform CD-UKF and CD-CKF. Also, we find that CD-EKF and CD-GBN-EKF achieve comparable accuracy as a result of similar working mechanisms. By Figure \ref{fig:Plot3}, we are able to observe the ARMSE for the filters, and see that CD-EKF and CD-GBN-EKF vastly outperform CD-UKF and CD-CKF, and their respective values are comparable. Further, CD-UKF and CD-CKF isn't able to complete a single run for larger values of $\delta \geq0.6$. This is because these filters cannot complete the integration of MDEs.

\subsubsection{Ill-conditioned Measurement Model}
For the ill-conditioned measurement model, we take the following measurement equation:
\begin{align}
z_k =
\begin{bmatrix}
1 & 1 \\
1 & 1+\sigma
\end{bmatrix}
x_k + v_k 
\end{align}
where $\sigma$ is a parameter of our choosing. $v_k$ is the same as before. For computational feasibility, we test ARMSE for CD-EKF and CD-GBN-EKF only. For larger values of $\sigma$, such as $\sigma = 10^{-2}$ and $\sigma = 10^{-4}$, both filters have comparable performance. Once $\sigma$ reaches $10^{-6}$, and becomes increasing ill-conditioned, we observe that CD-GBN-EKF outperforms CD-EKF. 

The numerical stability of the CD-GBN-EKF stems primarily from the absence of matrix inversion in its update mechanism. Conventional covariance-form filters, including the CD-EKF, require the inversion of the innovation covariance $S_k$, which becomes increasingly unstable as the measurement model increases ill-conditionedness:
\begin{align}
S_k &= H_k P_{k|k-1} H_k^\top + R_k, \\
K_k &= P_{k|k-1} H_k^\top S_k^{-1}. \label{eq:ekf_gain}
\end{align}
When $S_k$ is ill-conditioned, the inversion in~\eqref{eq:ekf_gain} can magnify numerical errors, degrade estimation accuracy, and even lead to loss of positive semi-definiteness in $P_{k|k}$. While square-root formulations alleviate some of these issues by propagating factors of the covariance matrix, they do not eliminate the dependence on $S_k^{-1}$ and therefore remain vulnerable in severely ill-conditioned settings.  

In contrast, the CD-GBN-EKF reformulates the update step using Gaussian Bayesian networks (GBNs). In a GBN, a multivariate Gaussian is decomposed into local conditional regressions of the form
\begin{align}
x_j &= \sum_{i=1}^{j-1} B_{ij} x_i + \epsilon_j, 
\qquad \epsilon_j \sim \mathcal{N}(0,v_j), \label{eq:gbn_regression}
\end{align}
where $B_{ij}$ are regression coefficients and $V_j$ are conditional variances. Measurement updates are performed through arc reversals and evidence entry, which only involve additive, multiplicative, and divisive operations on the parameters $(B,V)$. For example, when conditioning a child $x_j$ on its parent $x_i$, the updated variance follows
\begin{align}
v'_j = v_j + \frac{v_i}{B_{ij}^2}, \label{eq:gbn_update}
\end{align}
which requires no matrix inversion and is guaranteed to preserve positive semi-definiteness.  

From a numerical linear algebra perspective, this local representation distributes complexity across small-scale updates rather than relying on global inversion. Instead of computing $S_k^{-1}$ in one step, the CD-GBN-EKF incrementally conditions the joint distribution via~\eqref{eq:gbn_regression}--\eqref{eq:gbn_update}, thereby mitigating the amplification of roundoff errors and ensuring well-posed updates even in ill-conditioned regimes.  

As a result, the CD-GBN-EKF not only inherits the favorable stability of the CD-EKF in stiff stochastic systems, but also extends its robustness to scenarios where conventional filters fail due to ill-conditioned measurement models. This structural advantage makes Gaussian Bayesian network-based filtering a strong alternative for real-world applications where measurement matrices are nearly singular or poorly conditioned.

\begin{figure}
    \centering
    \includegraphics[width=0.42\textwidth]{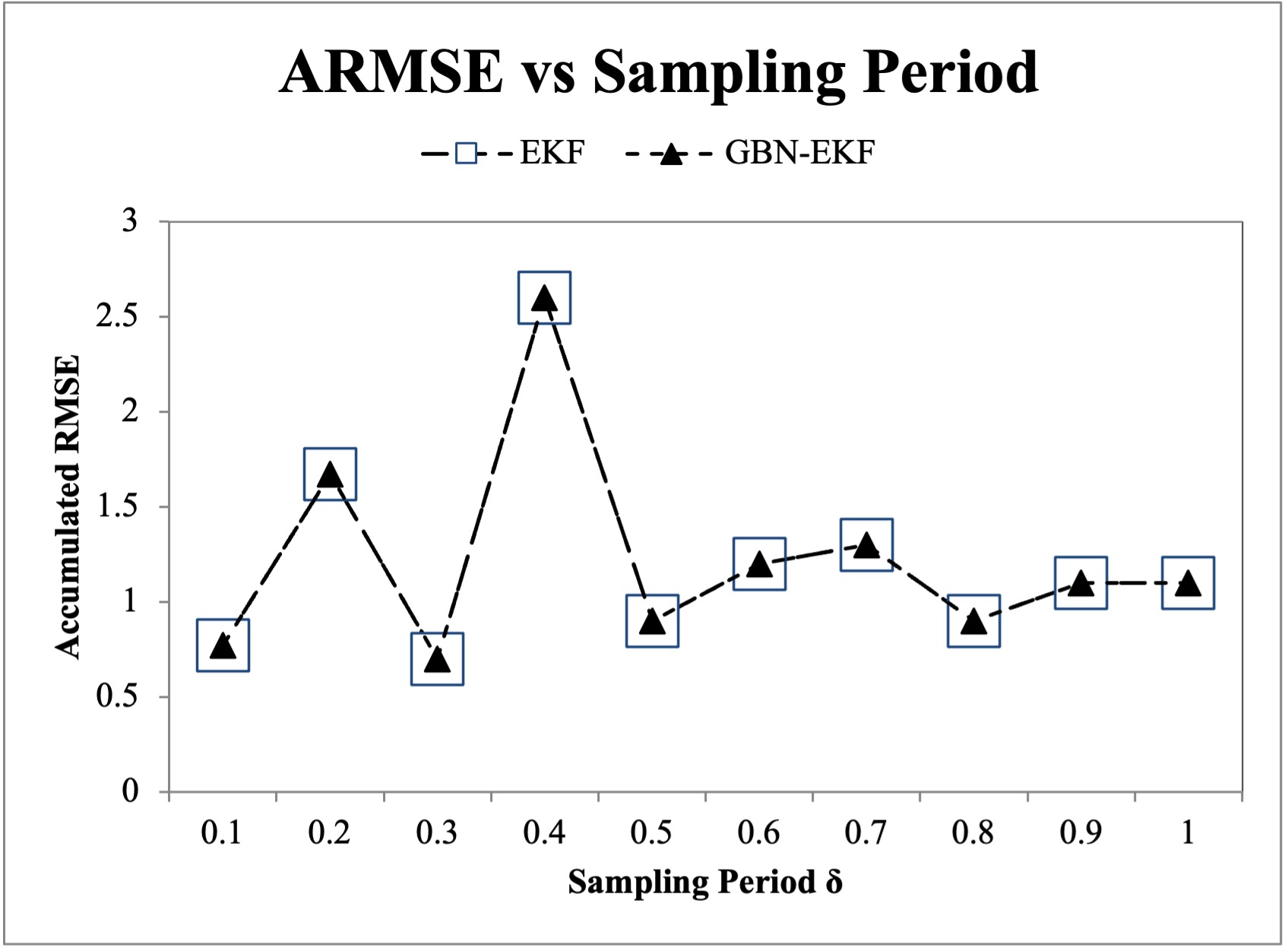}
    \caption{
        This figure plots ARMSE vs $\delta$ for the Van der Pol oscillator SDE with ill-conditioned measurements, $\sigma = 10^{-2}$
    }
    \label{fig:Plot4}
\end{figure}

\begin{figure}
    \centering
    \includegraphics[width=0.42\textwidth]{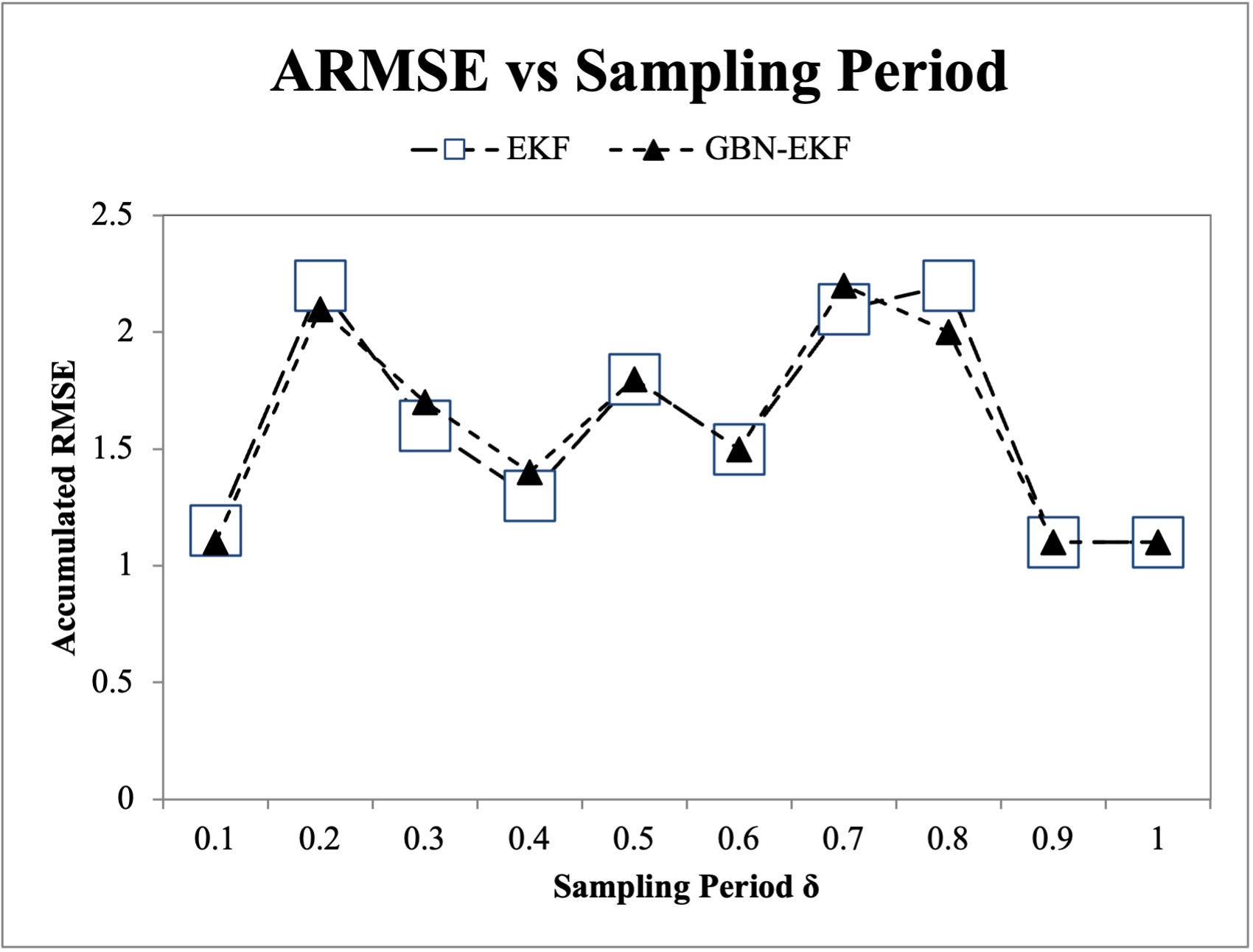}
    \caption{
        This figure plots ARMSE vs $\delta$ for the Van der Pol oscillator SDE with ill-conditioned measurements, $\sigma = 10^{-4}$
    }
    \label{fig:Plot5}
\end{figure}

\begin{figure}
    \centering
    \includegraphics[width=0.42\textwidth]{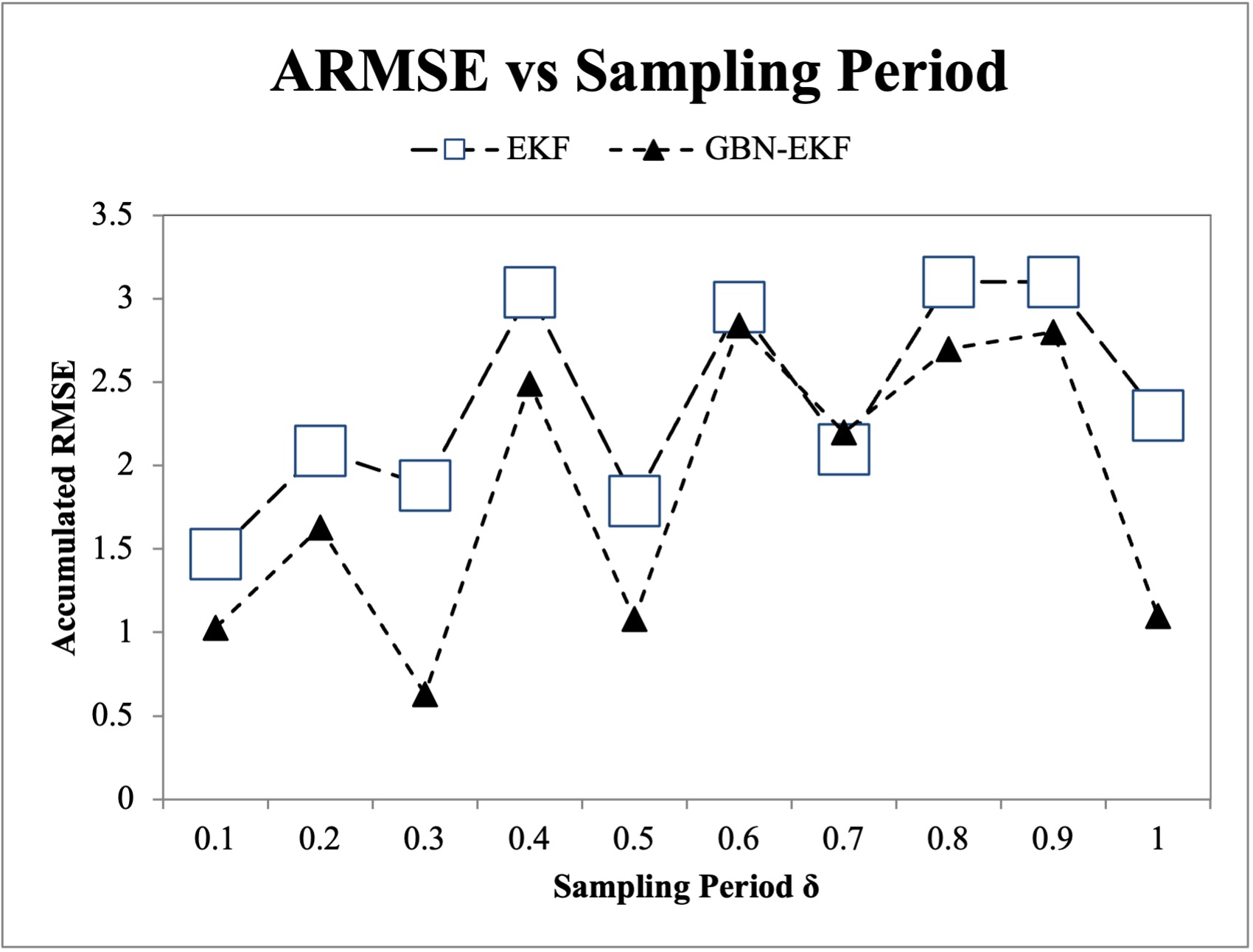}
    \caption{
        This figure plots ARMSE vs $\delta$ for the Van der Pol oscillator SDE with ill-conditioned measurements, $\sigma = 10^{-6}$
    }
    \label{fig:Plot6}
\end{figure}

\begin{figure}
    \centering
    \includegraphics[width=0.42\textwidth]{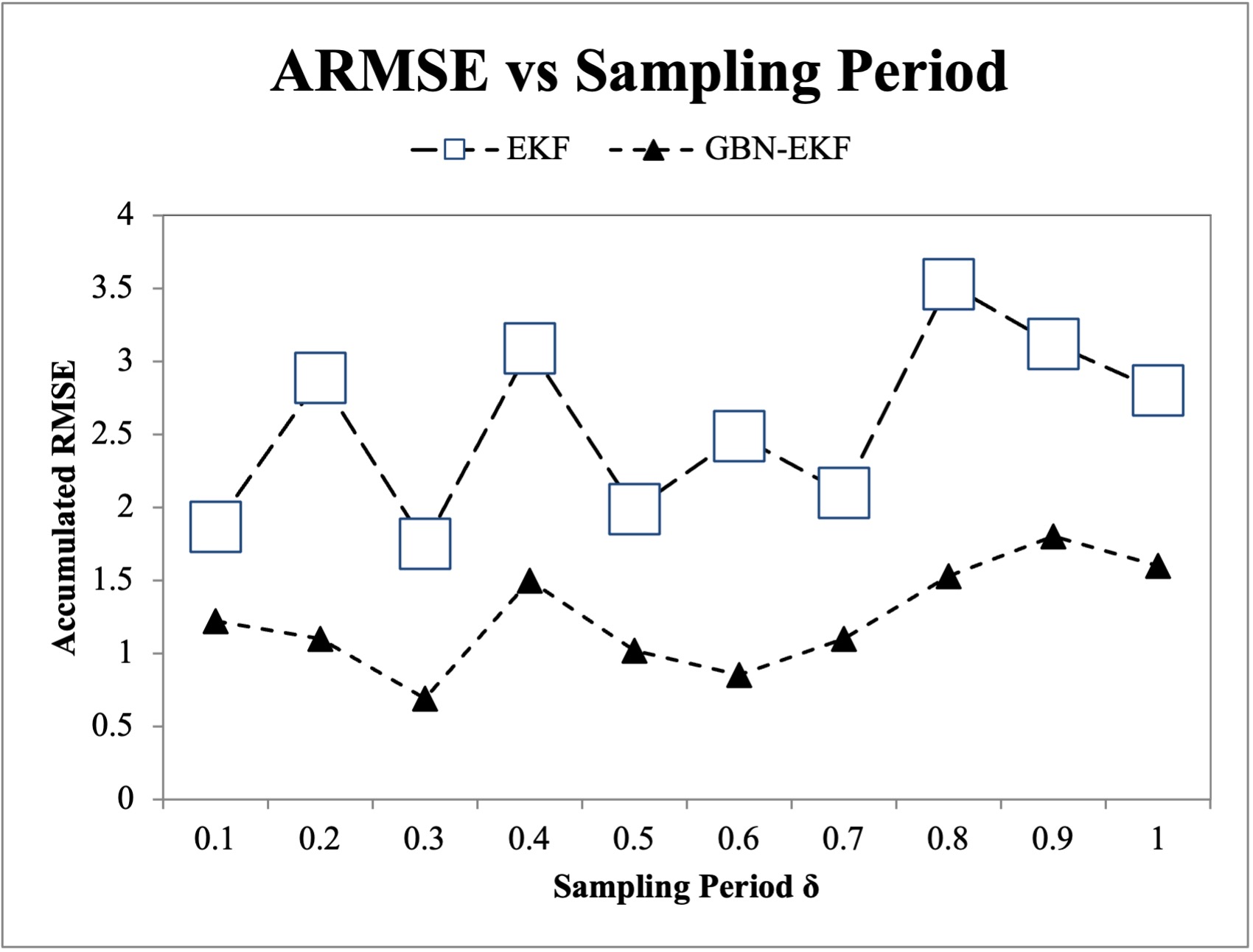}
    \caption{
        This figure plots ARMSE vs $\delta$ for the Van der Pol oscillator SDE with ill-conditioned measurements, $\sigma = 10^{-8}$
    }
    \label{fig:Plot7}
\end{figure}

\section{Conclusion}
This paper introduced a Gaussian Bayesian Network–based Extended Kalman Filter (CD-GBN-EKF) for continuous–discrete nonlinear state estimation in stiff systems, with a focus on ill-conditioned measurement models. By expressing the measurement update in a Gaussian Bayesian network and operating on local regression and conditional-variance parameters, the CD-GBN-EKF eliminates the explicit inversion of the innovation covariance. This structural change preserves positive semi-definiteness and improves numerical robustness when measurement matrices are poorly conditioned.

Numerical studies on the stiff Dahlquist and Van der Pol problems show that: (i) under well-conditioned measurements the CD-GBN-EKF matches the CD-EKF in accuracy, and (ii) under ill-conditioned measurements the CD-GBN-EKF consistently achieves lower ARMSE than the CD-EKF, while sigma-point methods (CD-UKF, CD-CKF) may degrade or fail due to instability of the underlying MDE propagation.

The proposed method inherits the first-order linearization of the EKF and may therefore be sensitive to severe nonlinearities or poor Jacobian quality. Future work can try to extend the CD-GBN-EKF through a square root formulation, and analyzing performance against existing square root methods in CD-EKF. 

\bibliography{references}

\end{document}